\documentclass{amsart}

\newtheorem{theorem}{Theorem}[section]
\newtheorem{lemma}[theorem]{Lemma}
\newtheorem{proposition}[theorem]{Proposition}
\newtheorem{corollary}[theorem]{Corollary}

\newtheorem{remark}[theorem]{Remark}

\theoremstyle{definition}
\newtheorem{definition}[theorem]{Definition}

\newcommand{\N}{\mathbb{N}}
\newcommand{\R}{\mathbb{R}}
\newcommand{\C}{\mathbb{C}}

\newcommand{\sig}{\sigma}
\newcommand{\Sig}{\Sigma}
\newcommand{\Degree}{d_{\sigma}}
\newcommand{\B}{\textnormal{B}_{\Sig}}
\newcommand{\K}{\textnormal{K}_{\Sig}}

\newcommand{\Lebesgue}{\mathcal{L}}

\begin{document}

\title[B-convex operator spaces]{B-convex operator spaces}

\author[Javier Parcet]{Javier Parcet}

\address{Department of Mathematics, Universidad Aut\'{o}noma
de Madrid.}

\email{javier.parcet@uam.es}


\maketitle

\begin{abstract}
The notion of $\mathrm{B}$-convexity for operator spaces, which a
priori depends on a set of parameters indexed by $\Sig$, is
defined. Some of the classical characterizations of this geometric
notion for Banach spaces are studied in this new context. For
instance, an operator space is $\B$-convex if and only if it has
$\Sig$-subtype. The class of uniformly non-$\Lebesgue^1(\Sig)$
operator spaces, which is also the class of $\B$-convex operator
spaces, is introduced. Moreover, an operator space having
non-trivial $\Sig$-type is $\B$-convex. However, the converse is
false. The row and column operator spaces are nice counterexamples
of this fact, since both are Hilbertian. In particular, this
result shows that a version of the Maurey-Pisier theorem does not
hold in our context. Some other examples of Hilbertian operator
spaces will be treated. In the last part of this paper, the
independence of $\B$-convexity with respect to $\Sig$ is studied.
This provides some interesting problems which will be posed.
\end{abstract}

\section{Introduction}
\label{Section1}

In the last thirty years, the notions of type and cotype of a
Banach space with respect to several orthonormal systems have been
deeply investigated. It turns out that, in order to study some
geometric properties of the Banach space considered, these notions
are very useful. In a series of papers with Jos\'{e} Garc\'{\i}a-Cuerva and
Jos\'{e} Manuel Marco, we have initiated a theory of type and cotype
for operator spaces. As it might be expected, the question is to
what extent these notions are related to the geometry of operator
spaces.

Our previous results can be summarized as follows. In \cite{GP1},
we study the validity of the Hausdorff-Young inequality for
vector-valued functions defined on a non-commutative compact
group. This inequality does not make sense when our functions take
values in a Banach space, we need an operator space structure on
it. This requirement goes back to Pisier's work \cite{P4}, where
non-commutative vector-valued Lebesgue spaces are studied. In
particular, this gives rise to the notions of Fourier type and
cotype of an operator space with respect to a non-commutative
compact group. We investigate in \cite{GP1} the basic properties
of this notion. The paper \cite{GMP} is devoted to the study of
the sharp Fourier type and cotype exponents of Lebesgue spaces and
Schatten classes. This is basic in the commutative theory.
However, the problem of finding the sharp exponents of a given
operator space is highly non-trivial, even for the simplest case
of Lebesgue spaces. When dealing with compact semisimple Lie
groups, we have solved part of this problem in \cite{GMP} using
the well-developed representation theory on this kind of groups.
Finally, the work \cite{GP2} deals with the general theory of type
and cotype for operator spaces. First, we define the notion of
quantized orthonormal system, which plays the role of the
orthonormal system in the classical setting. Then we introduce the
notions of type and cotype with respect to such a system. In
particular, this provides non-commutative extensions of the
notions of Rademacher, Steinhaus and Gaussian type and cotype. The
main result in \cite{GP2} is an operator space version of the
classical result of Kwapie\'n \cite{K}, which characterizes
Hilbert spaces up to isomorphism by means of vector-valued
orthogonal series. We give several approaches to this result in
the operator space setting, characterizing in such a way OH
operator spaces up to complete isomorphism by a type 2 / cotype 2
condition.

These previous results show that, also in the non-commutative
setting, there exists some interaction between the theory of type
and cotype and the geometry of operator spaces. Hence, it seems
that the next step in the process should be to study the notion of
B-convex operator space. Interest on B-convexity for Banach spaces
was generated in \cite{B}, where Beck studied certain strong law
of large numbers for Banach space valued random variables.
However, our motivation to study this notion for operator spaces
lies in some other characterizations of this geometric condition.
For instance, the following are equivalent for any Banach space
$B$: \\ \indent (a) $B$ is B-convex. \\ \indent (b) $B$ does not
contain $l_n^1$ uniformly. \\ \indent (c) $B$ has Rademacher
subtype. \\ \indent (d) $B$ has non-trivial Rademacher type. \\
\indent (e) $B$ is K-convex. \\ Giesy proved in \cite{G} the
equivalence between (a) and (b), while Pisier proved the others.
The proof of the equivalence between (a), (c) and (d) can be found
in \cite{P1}. Finally, the equivalence between B-convexity and
K-convexity was given in \cite{P2} and is much more complicated.
For some other characterizations of $\mathrm{B}$-convexity, the
reader is referred to \cite{PW}. Our aim is to study the validity
of these classical characterizations for operator spaces.

We shall follow a notation similar to that employed in \cite{GP2}.
Namely, let $(\Omega,\mathcal{M},\mu)$ be a probability measure
space with no atoms and let $\mathbf{d}_{\Sig} = \{\Degree: \sig
\in \Sig\}$ be a family of positive integers indexed by $\Sig$.
The \emph{quantized Steinhaus system of parameters} $(\Sig,
\mathbf{d}_{\Sig})$ is a collection $\mathcal{S}_{\Sig} =
\{\zeta^{\sig}: \Omega \rightarrow U(\Degree)\}_{\sig \in \Sig}$
of independent random unitary matrices, uniformly distributed on
the unitary group $U(\Degree)$ equipped with its normalized Haar
measure $\lambda_{\sig}$. Such a system lead us to the notion of
$\Sig$-type. There is a similar definition for the \emph{quantized
Rademacher system of parameters} $(\Sig, \mathbf{d}_{\Sig})$.
However, since we work with operator spaces (which are defined
over the complex field), it will be more convenient to deal with
the Steinhaus system. This does not make any difference since, as
it was mentioned in \cite{GP2}, Rademacher and Steinhaus
$\Sig$-type are equivalent notions.

The organization of the paper is as follows. In Section
\ref{Section2} we define the main notions we shall work with, such
as $\B$-convexity, $\Sig$-type, $\Sig$-subtype and uniformly
non-$\Lebesgue^1(\Sig)$ operator spaces. We also prove some basic
results which will be applied all throughout the paper. Some of
these notions are equivalent. Namely, an operator space has
$\Sig$-subtype if and only if it is uniformly
non-$\Lebesgue^1(\Sig)$ and each of these properties are
equivalent to the condition of being $\B$-convex. The proof of
this is the content of Section \ref{Section3}. In Section
\ref{Section4} we work with certain tensor submultiplicativity
which generalizes the classical submultiplicativity of the
constants involved in the theory, see \cite{P1} for the details.
This will be useful in order to see the validity in our context,
for the main sets of parameters $(\Sig, \mathbf{d}_{\Sig}) $, of
some well-known cla\-ssical results. In Section \ref{Section5} we
show that having non-trivial $\Sig$-type is stronger than
$\B$-convexity. Moreover, we provide examples of $\B$-convex
Hilbertian operator spaces failing the non-trivial $\Sig$-type
condition. Some other interesting examples will be treated.
Finally, Section \ref{Section6} is devoted to the study of the
dependence of $\B$-convexity with respect to $(\Sig,
\mathbf{d}_{\Sig})$. First, we prove that there is no dependence
when we work only with those sets of parameters having
$\mathbf{d}_{\Sig}$ unbounded. Then, we introduce the notion of
$\K$-convexity, and we show that the independence with respect to
$(\Sig, \mathbf{d}_{\Sig})$ is equivalent to see that an operator
space is $\B$-convex if and only if it is $\K$-convex. Also a
sufficient condition for the independence, in terms of the way in
which $S_n^1$ embeds in $S^2(l^1)$, is given.

All throughout this paper, some basic notions of operator spaces
and vector-valued Schatten classes will be assumed. The main
results that we will be using can be found in \cite{P4}. Also, the
main results in the theory of type for operator spaces will be
assumed. The reader is referred to \cite{GP2} for a brief summary
of them. Finally, along this work we make a slight abuse of
notation since sometimes, when there is no risk of confusion, we
shall write $\Sig$ to denote the set of parameters $(\Sig,
\mathbf{d}_{\Sig})$.

\section{The main definitions}
\label{Section2}

Given an operator space $E$, a set of parameters $(\Sig,
\mathbf{d}_{\Sig})$ and $1 \le p < \infty$, we define the spaces
$\Lebesgue_E^p(\Sig)$ as follows
\begin{eqnarray*}
\Lebesgue_E^p(\Sig) & = & \Big\{ A \in \prod_{\sig \in \Sig}
M_{\Degree} \otimes E: \,\ \|A\|_{\Lebesgue_E^p(\Sig)} = \Big(
\sum_{\sig \in \Sig} \Degree \|A^{\sig}\|_{S_{\Degree}^p(E)}^p
\Big)^{1/p} < \infty \Big\} \\ \Lebesgue_E^{\infty}(\Sig) & = &
\Big\{ A \in \prod_{\sig \in \Sig} M_{\Degree} \otimes E: \,\
\|A\|_{\Lebesgue_E^{\infty}(\Sig)} = \sup_{\sig \in \Sig}
\|A^{\sig}\|_{S_{\Degree}^{\infty}(E)} < \infty \Big\},
\end{eqnarray*}
where we write $M_n$ for the vector space of $n \times n$ complex
matrices and $S_n^p(E)$ stands for the Schatten $p$-class over
$M_n$ with values in $E$. We shall also use the
infinite-dimensional Schatten classes $S^p(E)$ with values in $E$,
the reader is referred to \cite{P4} for a precise definition of
these spaces. $\Lebesgue^p(\Sig)$ will denote the case $E = \C$.
We impose on $\Lebesgue_E^p(\Sig)$ its natural operator space
structure, see \cite{GP1} or Chapter 2 of \cite{P4} for the
details. Now let $1 \le p \le \infty$ and let $\Gamma$ be any
subset of $\Sig$, if $\zeta^{\sig}$ are the elements of the
quantized Steinhaus system (as defined in the introduction), then
we define the mapping $\mathrm{T}_p(\Gamma, E)$ by the relation
$$A \in \Lebesgue_E^p(\Gamma) \longmapsto \sum_{\sig \in \Gamma}
\Degree \mbox{tr}(A^{\sig} \zeta^{\sig}) \in L_E^2(\Omega).$$ We
shall denote $\mathrm{T}_p(\Sig,E)$ by $\mathrm{T}_p(E)$. For
$\Gamma$ finite, we also define the number $$\displaystyle
\Delta_{\Gamma} = \sum_{\sig \in \Gamma} \Degree^2.$$

\begin{definition} \label{Definition-Type}
Let $1 \le p \le 2$, we say that an operator space $E$ has
$\Sig$-\emph{type} $p$ if the mapping $\mathrm{T}_p(E)$ is
completely bounded.
\end{definition}

\begin{remark} \label{Relevant-Parameters}
\emph{This notion of type depends on the election of the quantized
Steinhaus system $\mathcal{S}_{\Sig}$ we are working with. If we
take $\Sig_0 = \N$ and $\Degree = 1$ for all $\sig \in \Sig_0$,
then we work with the classical Steinhaus system and Definition
\ref{Definition-Type} provides a completely bounded version of the
classical definition of type. We shall refer to these parameters
as the \emph{classical set of parameters} $\Sig_0$. One could
think that this is the right definition of type in the operator
space setting and that there is no reason to introduce all the
quantized Steinhaus systems in the theory. But in fact, those
quantized Steinhaus systems $\mathcal{S}_{\Sig}$ with
$\mathbf{d}_{\Sig}$ unbounded are sometimes the right systems to
work with. One example of this assertion is given in \cite{GP2},
where we show that for those systems we can obtain an operator
space version of Kwapie\'n's theorem with weaker hypothesis than
for the classical Steinhaus system. Another example will be given
in Section \ref{Section6}, where we shall prove that the notion of
$\B$-convexity is independent of $\Sig$ whenever
$\mathbf{d}_{\Sig}$ is unbounded.}
\end{remark}

\begin{remark} \label{Khintchine-Kahane}
\emph{Let $\mathcal{S}_{\Sig}^q(E)$ be the closure in
$L_E^q(\Omega)$ of the subspace spanned by the entries
$\{\zeta_{ij}^{\sig}: \sig \in \Sig, \,\ 1 \le i,j \le \Degree\}$
of the functions of $\mathcal{S}_{\Sig}$ with $E$-valued
coefficients. A version of the Khintchine-Kahane inequalities for
random matrices given in \cite{MP} implies that the norm of
$\mathcal{S}_{\Sig}^{q_1}(E)$, regarded as a Banach space, is
equivalent to that of $\mathcal{S}_{\Sig}^{q_2}(E)$ whenever $1
\le q_1, q_2 < \infty$. In particular, given $1 \le p \le 2$ and
$1 \le q < \infty$, the validity of the inequality
\begin{equation} \label{Independence}
\Big\| \sum_{\sig \in \Sig} \Degree \mbox{tr}(A^{\sig}
\zeta^{\sig}) \Big\|_{L_E^q(\Omega)} \le c \,\ \Big( \sum_{\sig
\in \Sig} \Degree \|A^{\sig}\|_{S_{\Degree}^p(E)}^p \Big)^{1/p}
\end{equation}
does not depend on the value of $q$. However, for $1 \le q_1 \neq
q_2 < \infty$, Pisier showed in \cite{P4} that in general
$\mathcal{S}_{\Sig}^{q_1}(E)$ and $\mathcal{S}_{\Sig}^{q_2}(E)$
are not completely isomorphic as operator spaces with their
natural operator space structure. Therefore, in contrast with
(\ref{Independence}), each election of $1 \le q < \infty$ in
Definition \ref{Definition-Type} gives different notions of
$\Sig$-type.}
\end{remark}

\begin{remark} \label{Independence-Type}
\emph{The absence of Khintchine inequalities for operator spaces
forces us to choose an exponent $q$ in the definition of
$\Sig$-type. Our election differs from that of \cite{GP2}. The
reason is that in \cite{GP2} (where we took $q = p'$ in the
definition of $\Sig$-type $p$) the aim was to have a unified
theory of type and cotype for uniformly bounded quantized
orthonormal systems, while here the election $q = 2$ facilitates
our work. However, although the notion of $\Sig$-type given here
is not the same as in \cite{GP2}, we shall mainly be concerned
with the notion of \emph{non-trivial $\Sig$-type} for which the
election of $q$ does not matter! Namely, let $1 \le q_1 < q_2 <
\infty$ and $1 < p \le 2$. Let us consider the mapping
$\mathrm{T}_p^q(E)$ given by $$A \in \Lebesgue_E^p(\Sig)
\longmapsto \sum_{\sig \in \Sig} \Degree \mbox{tr}(A^{\sig}
\zeta^{\sig}) \in L_E^q(\Omega).$$ The complete boundedness of
$\mathrm{T}_p^{q_2}(E)$ obviously implies the complete boundedness
of $\mathrm{T}_p^{q_1}(E)$. Conversely, let us assume that
$\mathrm{T}_p^{q_1}(E)$ is completely bounded. Then, since the
$cb$-norm of $\mathrm{T}_1^{\infty}(E)$ is 1 for any operator
space $E$, by complex interpolation there exists some $1 < r < p$
such that $\mathrm{T}_r^{q_2}(E)$ is also completely bounded. That
is, the notion of non-trivial $\Sig$-type does not depend on the
election of $q$ in the definition of $\Sig$-type.}
\end{remark}

We also generalize to our context some other notions of the
commutative theory. For that we shall need the numbers
$\Delta_{\Gamma}$ defined above for each finite $\Gamma \subset
\Sig$. Let $E$ be an operator space and let $(\Sig,
\mathbf{d}_{\Sig})$ be any set of parameters. We say that $E$ has
\emph{$\Sig$-subtype} if there exists a finite subset $\Gamma$ of
$\Sig$, such that $$\|\mathrm{T}_2(\Gamma, E)\|_{cb} <
\Delta_{\Gamma}^{1/2}.$$

Now we define $\B$-convex and uniformly non-$\Lebesgue^1(\Sig)$
operator spaces. As in the previous definitions, when dealing with
the classical set of parameters $\Sig_0$, we obtain a completely
bounded version of the classical notion.

\begin{definition} \label{Definitions}
Let $E$ be an operator space and let us fix $(\Sig,
\mathbf{d}_{\Sig})$.
\begin{itemize}
\item $E$ is said to be \emph{$\B$-convex} if there exists a
finite subset $\Gamma$ of $\Sig$ and $0 < \delta \le 1$ such that,
for any family $\{A^{\sig} \in M_{\Degree} \otimes S^2(E)\}_{\sig
\in \Gamma}$, we have $$\frac{1}{\Delta_{\Gamma}} \inf_{B^{\sig}
unitary} \Big\|\sum_{\sig \in \Gamma} \Degree \mbox{tr}(A^{\sig}
B^{\sig}) \Big\|_{S^2(E)} \le (1 - \delta) \,\ \max_{\sig \in
\Gamma}\|A^{\sig}\|_{S_{\Degree}^{\infty}(S^2(E))}.$$
\item $E$ is said to \emph{contain $\Lebesgue^1(\Gamma)$
$\lambda$-uniformly} if, for each finite subset $\Gamma$ of
$\Sig$, there exists a subspace $F_{\Gamma}$ of $S^2(E)$ and a
linear isomorphism $\Lambda_{\Gamma}: \Lebesgue^1(\Gamma)
\rightarrow F_{\Gamma}$ such that $$\|\Lambda_{\Gamma}\|_{cb}
\|\Lambda_{\Gamma}^{-1}\| \le \lambda.$$ $E$ is called a
\emph{uniformly non-$\Lebesgue^1(\Sig)$} operator space if it does
not contain $\Lebesgue^1(\Gamma)$ $\lambda$-uniformly for some
$\lambda > 1$.
\end{itemize}
\end{definition}

\begin{remark}
\emph{The reader could expect that, in the given definition of
containing $\Lebesgue^1(\Gamma)$ $\lambda$-uniformly, we should
require that $\|\Lambda_{\Gamma}\|_{cb}
\|\Lambda_{\Gamma}^{-1}\|_{cb} \le \lambda$. However we have given
an intermediate notion between that condition and the classical
notion, which uses the Banach-Mazur distance. The reason for that
election will become clear in Theorem \ref{Equivalences}.}
\end{remark}

\begin{remark} \label{Uniformly-Classical}
\emph{For the classical set of parameters $\Sig_0$, the given
definition of containing $\Lebesgue^1(\Gamma)$ $\lambda$-uniformly
can be rephrased by saying that the space $S^2(E)$ contains
$l_n^1$ $\lambda$-uniformly in the Banach space sense. Namely, the
$cb$-norm of $\Lambda_{\Gamma}$ coincides  in this case with its
operator norm since it is defined on a $\max$ operator space.}
\end{remark}

\begin{remark} \label{Remark-p-Independence}
\emph{We shall see in Section \ref{Section3} that Definition
\ref{Definitions} does not change if we put $S^p(E)$ instead of
$S^2(E)$ for any $1 < p < \infty$.}
\end{remark}

We need to prove the following result which will be used sometimes
in this paper and for which we have not found any reference.

\begin{lemma} \label{cb-distance}
Let $E$ be an operator space and let $n$ be a positive integer.
\begin{itemize}
\item[$(\mathrm{a})$] If $1 \le p < q \le \infty$, then $\|A\|_{S_n^p(E)} \le n^{1/p -
1/q} \|A\|_{S_n^q(E)}$.
\item[$(\mathrm{b})$] If $\|A\|_{S_n^1(E)} = \sqrt{n} \,\ \|A\|_{S_n^2(E)}$, then
$\|A\|_{S_n^1(E)} = n \,\ \|A\|_{S_n^{\infty}(E)}$.
\end{itemize}
\end{lemma}

\begin{proof}
If $p = 1$ and $q = \infty$, then (a) follows easily from
Corollary 9.8 of \cite{P3}. The general case follows by complex
interpolation. Let us prove (b). By homogeneity, we assume without
lost of generality that $\|A\|_{S_n^2(E)} = 1$. On the other hand,
by Theorem 1.5 of \cite{P4}, we know that $$\|A\|_{S_n^2(E)} =
\inf_{A = \alpha B \beta} \|\alpha\|_{S_n^{4}}
\|B\|_{S_n^{\infty}(E)} \|\beta\|_{S_n^4} = 1$$ where $\alpha,
\beta \in M_n$ and $B \in M_n \otimes E$. Moreover, if $F =
\mbox{span} \{A_{ij}: 1 \le i,j \le n\}$, we can take $B \in M_n
\otimes F$. In particular, for all $k \ge 1$ there exist
$\alpha_k, \beta_k \in M_n$ and $B_k \in M_n \otimes F$ such that
$A = \alpha_k B_k \beta_k$, $1 \le \|\alpha_k\|_{S_n^4} < 1 + 1/k$
and $\|B_k\|_{S_n^{\infty}(F)} = \|\beta_k\|_{S_n^4} = 1$. By the
finite-dimensionality of $F$, we know that the sequence
$(\alpha_k, B_k, \beta_k)$ belongs to a compact subset of $S_n^4
\times S_n^{\infty}(E) \times S_n^4$. Thus, there exist $\alpha_0,
\beta_0 \in M_n$ and $B_0 \in M_n \otimes F$ such that $A =
\alpha_0 B_0 \beta_0$ and $\|\alpha_0\|_{S_n^4} =
\|B_0\|_{S_n^{\infty}(F)} = \|\beta_0\|_{S_n^4} =1$. But then,
again by Theorem 1.5 of \cite{P4}, $$\|\alpha_0\|_{S_n^2}
\|\beta_0\|_{S_n^2} \ge \|A\|_{S_n^1(E)} = \sqrt{n}.$$ Moreover,
taking $p=2$ and $q=4$, (a) gives $\|\alpha_0\|_{S_n^2},
\|\beta_0\|_{S_n^2} \le n^{1/4}$. In summary, we get
$\|\alpha_0\|_{S_n^2} = \|\beta_0\|_{S_n^2} = n^{1/4}$ and
$\|\alpha_0\|_{S_n^4} = \|\beta_0\|_{S_n^4} = 1$. Then, it is
well-known that there exists $U,V \in U(n)$ such that $\alpha_0 =
n^{-1/4} U$ and $\beta_0 = n^{-1/4} V$. Therefore,
$\|A\|_{S_n^{\infty}(E)} \le \|\alpha_0\|_{S_n^{\infty}}
\|B_0\|_{S_n^{\infty}(E)} \|\beta_0\|_{S_n^{\infty}} = n^{-1/2}$.
This gives $\|A\|_{S_n^1(E)} \ge n \,\ \|A\|_{S_n^{\infty}(E)}$,
the reverse inequality follows from (a). This completes the proof.
\end{proof}

\section{The equivalent notions}
\label{Section3}

In this section we show the equivalence between some of the
notions previously defined. For that purpose we begin by fixing a
set of parameters $(\Sig, \mathbf{d}_{\Sig})$. First, we need to
prove a technical lemma.

\begin{lemma} \label{Ultraproducts}
Let $F$ be an operator space and let $\Gamma$ be a finite subset
of $\Sig$. Let us suppose that, for all $\varepsilon > 0$, there
exist a family of matrices $X_{\varepsilon}^{\sig} \in M_{\Degree}
\otimes F$ with $\sig \in \Gamma$ and such that
\begin{itemize}
\item[$1.$] $\displaystyle \Big( \int_{\Omega} \Big\| \sum_{\sig \in
\Gamma} \Degree \textnormal{tr}(X_{\varepsilon}^{\sig}
\zeta^{\sig}(\omega)) \Big\|_{F}^2 d \mu(\omega) \Big)^{1/2} \ge
\Delta_{\Gamma} - \varepsilon$.
\item[$2.$] $\displaystyle \Big( \sum_{\sig \in \Gamma} \Degree
\|X_{\varepsilon}^{\sig}\|_{S_{\Degree}^2(F)}^2 \Big)^{1/2} =
\Delta_{\Gamma}^{1/2}$.
\end{itemize}
Then $\displaystyle \max_{\sig \in \Gamma}
\|X_{\varepsilon}^{\sig}\|_{S_{\Degree}^{\infty}(F)} \le 1 +
\xi(\varepsilon)$, with $\xi(\varepsilon) \rightarrow 0^+$ as
$\varepsilon \rightarrow 0^+$.
\end{lemma}

\begin{proof}
Let $\mathcal{U}$ be an ultrafilter on $\R_+$ containing all the
intervals $(0, \varepsilon)$ with $\varepsilon > 0$ and let
$F_{\mathcal{U}}$ be the corresponding ultraproduct operator
space. Then we define $X^{\sig} =
(X_{\varepsilon}^{\sig})_{\mathcal{U}}$ for $\sig \in \Gamma$.
That is, $X^{\sig} \in M_{\Degree} \otimes F_{\mathcal{U}}$. We
obviously have $$\Big( \int_{\Omega} \Big\| \sum_{\sig \in \Gamma}
\Degree \textnormal{tr}(X^{\sig} \zeta^{\sig}(\omega))
\Big\|_{F_{\mathcal{U}}}^2 d \mu(\omega) \Big)^{1/2} \ge
\Delta_{\Gamma}$$ and also $$\Big( \sum_{\sig \in \Gamma} \Degree
\|X^{\sig}\|_{S_{\Degree}^2(F_{\mathcal{U}})}^2 \Big)^{1/2} =
\Delta_{\Gamma}^{1/2}.$$ By H\"{o}lder inequality and Lemma
\ref{cb-distance}, we can write
\begin{eqnarray*}
\Delta_{\Gamma}^2 & \le & \int_{\Omega} \Big( \sum_{\sig \in
\Gamma} \Degree \|X^{\sig}\|_{S_{\Degree}^1(F_{\mathcal{U}})}
\|\zeta^{\sig}(\omega)\|_{S_{\Degree}^{\infty}} \Big)^2 d
\mu(\omega) \\ & = & \Big( \sum_{\sig \in \Gamma} \Degree
\|X^{\sig}\|_{S_{\Degree}^1(F_{\mathcal{U}})} \Big)^2 \le \Big(
\sum_{\sig \in \Gamma} \Degree \sqrt{\Degree}
\|X^{\sig}\|_{S_{\Degree}^2(F_{\mathcal{U}})} \Big)^2 \\ & \le &
\Delta_{\Gamma} \sum_{\sig \in \Gamma} \Degree
\|X^{\sig}\|_{S_{\Degree}^2(F_{\mathcal{U}})}^2 =
\Delta_{\Gamma}^2.
\end{eqnarray*}
In particular, in this case, Lemma \ref{cb-distance} and H\"{o}lder
inequality are equalities. Therefore, we obtain
$$\|X^{\sig}\|_{S_{\Degree}^1(F_{\mathcal{U}})} = \sqrt{\Degree}
\|X^{\sig}\|_{S_{\Degree}^2(F_{\mathcal{U}})} \qquad \mbox{and}
\qquad \Degree = c_0 \sqrt{\Degree}
\|X^{\sig}\|_{S_{\Degree}^2(F_{\mathcal{U}})}$$ for some positive
constant $c_0$ and any $\sig \in \Gamma$. Now, Lemma
\ref{cb-distance} gives
$\|X^{\sig}\|_{S_{\Degree}^{\infty}(F_{\mathcal{U}})} = 1 / c_0$.
But $c_0 = 1$ since $$\frac{1}{c_0^2} \,\ \Delta_{\Gamma} =
\sum_{\sig \in \Gamma} \Degree
\|X^{\sig}\|_{S_{\Degree}^2(F_{\mathcal{U}})}^2 =
\Delta_{\Gamma}.$$ So, we have $\max_{\sig \in \Gamma}
\|X^{\sig}\|_{S_{\Degree}^{\infty}(F_{\mathcal{U}})} = 1$.
Finally, by the isometry $S_{\Degree}^{\infty}(F_{\mathcal{U}}) =
S_{\Degree}^{\infty}(F)_{\mathcal{U}}$ (see Chapter 5 of
\cite{P4}) and the definition of ultraproduct operator space, the
result follows. This completes the proof.
\end{proof}

\begin{theorem} \label{Equivalences}
Let $E$ be an operator space and let us fix a set of parameters
$(\Sig, \mathbf{d}_{\Sig})$. Then, the following are equivalent:
\begin{itemize}
\item[$(\mathrm{a})$] $E$ has $\Sig$-subtype.
\item[$(\mathrm{b})$] $E$ is $\B$-convex.
\item[$(\mathrm{c})$] $E$ is uniformly non-$\Lebesgue^1(\Sig)$.
\end{itemize}
\end{theorem}

\begin{proof}
Let us suppose that $E$ has $\Sig$-subtype, we shall see that $E$
is $\B$-convex. We know that $\|\mathrm{T}_2(\Gamma, E)\|_{cb} =
(1 - \delta) \,\ \Delta_{\Gamma}^{1/2}$ for some $\Gamma \subset
\Sig$ finite and some $0 < \delta \le 1$. Hence, we can write
\begin{eqnarray*}
\frac{1}{\Delta_{\Gamma}} \inf_{B^{\sig} unitary} \Big\|
\sum_{\sig \in \Gamma} \Degree \mbox{tr}(A^{\sig} B^{\sig})
\Big\|_{S^2(E)} & \le & \frac{1}{\Delta_{\Gamma}} \Big\|
\sum_{\sig \in \Gamma} \Degree \mbox{tr}(A^{\sig} \zeta^{\sig})
\Big\|_{S^2(L_E^2(\Omega))} \\ & \le & \frac{(1 -
\delta)}{\Delta_{\Gamma}^{1/2}} \,\
\|A\|_{S^2(\Lebesgue_E^2(\Gamma))}.
\end{eqnarray*}
Then the result follows since, by Lemma \ref{cb-distance}, we have
$$\|A\|_{\Lebesgue_{S^2(E)}^2(\Gamma)} \le \Delta_{\Gamma}^{1/2}
\,\ \max_{\sig \in \Gamma}
\|A^{\sig}\|_{S_{\Degree}^{\infty}(S^2(E))}.$$

Now, to see that $\B$-convex operator spaces are uniformly
non-$\Lebesgue^1(\Sig)$, we assume that $E$ contains
$\Lebesgue^1(\Gamma)$ $\lambda$-uniformly for all $\lambda > 1$
and we have to see that $E$ is not $\B$-convex. We know that, for
all $\lambda > 1$ and all $\Gamma \subset \Sig$ finite, there
exists a subspace $F_{\Gamma}$ of $S^2(E)$ and some isomorphism
$\Lambda_{\Gamma}: \Lebesgue^1(\Gamma) \rightarrow F_{\Gamma}$ $$A
\in \Lebesgue^1(\Gamma) \longmapsto \sum_{\sig \in \Gamma}
\sum_{i,j =1}^{\Degree} a_{ij}^{\sig} x_{ij}^{\sig} \in F_{\Gamma}
\qquad \mbox{where} \ \ A^{\sig} = (a_{ij}^{\sig}),$$ such that
$\|\Lambda_{\Gamma}\|_{cb} = 1$ and $\|\Lambda_{\Gamma}^{-1}\| \le
\lambda$. On the other hand, if $\sig \in \Gamma$ and we define
the matrix $X^{\sig} = \Degree^{-1} (x_{ij}^{\sig})$, we have
\begin{equation} \label{Norm-Estimate}
\|X^{\sig}\|_{S_{\Degree}^{\infty}(S^2(E))} = \|\mbox{tr}(X^{\sig}
\cdot)\|_{cb(S_{\Degree}^1, S^2(E))} \le \|\Lambda_{\Gamma}\|_{cb}
= 1.
\end{equation}
Hence, by the estimate for the norm of $\Lambda_{\Gamma}^{-1}$ and
(\ref{Norm-Estimate}), we obtain $$\frac{1}{\Delta_{\Gamma}}
\inf_{B^{\sig} unitary} \Big\| \sum_{\sig \in \Gamma} \Degree
\mbox{tr}(X^{\sig} B^{\sig}) \Big\|_{S^2(E)} \ge \frac{1}{\lambda}
\max_{\sig \in \Gamma}
\|X^{\sig}\|_{S_{\Degree}^{\infty}(S^2(E))}$$ since
$\|B\|_{\Lebesgue^1(\Gamma)} = \Delta_{\Gamma}$ whenever $B^{\sig}
\in U(\Degree)$ for all $\sig \in \Gamma$. In particular, taking
$\lambda \rightarrow 1^+$, we conclude that the operator space $E$
is not $\mbox{B}_{\Sig}$-convex.

Finally, let us assume that $\|\mathrm{T}_2(\Gamma, E)\|_{cb} =
\Delta_{\Gamma}^{1/2}$ for all $\Gamma \subset \Sig$ finite. We
have to see that $E$ contains $\Lebesgue^1(\Gamma)$
$\lambda$-uniformly for all $\lambda > 1$. By Lemma 1.7 of
\cite{P4}, we know that $\|\mathrm{T}_2(\Gamma, E)\|_{cb} =
\|\mathrm{T}_2(\Gamma, S^2(E))\|$. In particular, for all
$\varepsilon > 0$ there exists a family of matrices
$X_{\varepsilon}^{\sig} \in M_{\Degree} \otimes S^2(E)$ such that
$$\int_{\Omega} \Big\| \sum_{\sig \in \Gamma} \Degree
\mbox{tr}(X_{\varepsilon}^{\sig} \zeta^{\sig}(\omega))
\Big\|_{S^2(E)}^2 d \mu(\omega) \ge \Delta_{\Gamma}^2 -
\varepsilon$$ and such that $$\sum_{\sig \in \Gamma} \Degree
\|X_{\varepsilon}^{\sig}\|_{S_{\Degree}^2(S^2(E))}^2 =
\Delta_{\Gamma}.$$ Moreover, by H\"{o}lder inequality we obtain
$$\Big\| \sum_{\sig \in \Gamma} \Degree
\mbox{tr}(X_{\varepsilon}^{\sig} \zeta^{\sig} (\omega))
\Big\|_{S^2(E)} \le \sum_{\sig \in \Gamma} \Degree
\|X_{\varepsilon}^{\sig}\|_{S_{\Degree}^2(S^2(E))}
\|\zeta^{\sig}(\omega)\|_{S_{\Degree}^2} \le \Delta_{\Gamma}.$$
That is, if we set $\displaystyle f_{\varepsilon}(\omega) = \Big\|
\sum_{\sig \in \Gamma} \Degree \mbox{tr}(X_{\varepsilon}^{\sig}
\zeta^{\sig} (\omega)) \Big\|_{S^2(E)}^2$ for $\omega \in \Omega$,
then we have $$0 \le f_{\varepsilon} \le \Delta_{\Gamma}^2 \qquad
\mbox{and} \qquad \Delta_{\Gamma}^2 - \varepsilon \le
\int_{\Omega} f_{\varepsilon}(\omega) d \mu(\omega) \le
\Delta_{\Gamma}^2.$$ In particular, $\mu\{\omega \in \Omega:
f_{\varepsilon}(\omega) < \Delta_{\Gamma}^2 - k \varepsilon\} \le
1/k$ for all $k \ge 1$. On the other hand, if we fix $U_0^{\sig}
\in U(\Degree)$ for any $\sig \in \Gamma$, we define
$$\mathcal{U}_0(\sig, \delta) = \{U^{\sig} \in U(\Degree): \,\
\|U^{\sig} - U_0^{\sig}\|_{S_{\Degree}^2} < \delta\}.$$ Then we
recall that, by the independence of the random matrices
$\zeta^{\sig}$ and their uniform distribution in $U(\Degree)$ with
respect to the normalized Haar measure $\lambda_{\sig}$ in
$U(\Degree)$, we have $$\mu\{\omega \in \Omega: \,\
\zeta^{\sig}(\omega) \in \mathcal{U}_0(\sig, \delta), \,\ \sig \in
\Gamma\} = \prod_{\sig \in \Gamma}
\lambda_{\sig}(\mathcal{U}_0(\sig, \delta)) > 0.$$ Therefore, by
choosing $k_0(\delta)$ such that $k_0(\delta)^{-1} < \mu\{\omega
\in \Omega: \,\ \zeta^{\sig}(\omega) \in \mathcal{U}_0(\sig,
\delta), \,\ \sig \in \Gamma\}$, we obtain the following
inequality $$\mu \{\omega \in \Omega: f_{\varepsilon}(\omega) <
\Delta_{\Gamma}^2 - k_0(\delta) \varepsilon\} < \mu\{\omega \in
\Omega: \,\ \zeta^{\sig}(\omega) \in \mathcal{U}_0(\sig, \delta),
\,\ \sig \in \Gamma\}.$$ That is, there exists some $\omega_0 \in
\Omega$ such that $\zeta^{\sig}(\omega_0) \in \mathcal{U}_0(\sig,
\delta)$ for all $\sig \in \Gamma$ and such that
$f_{\varepsilon}(\omega_0) \ge \Delta_{\Gamma}^2 - k_0(\delta)
\varepsilon$. These two properties give us the following sequence
of inequalities
\begin{eqnarray*}
\lefteqn{\sqrt{\Delta_{\Gamma}^2 - k_0(\delta) \varepsilon}} \\ &
\le & \Big\| \sum_{\sig \in \Gamma} \Degree
\mbox{tr}(X_{\varepsilon}^{\sig} [\zeta^{\sig}(\omega_0) -
U_0^{\sig}]) \Big\|_{S^2(E)} + \Big\| \sum_{\sig \in \Gamma}
\Degree \mbox{tr}(X_{\varepsilon}^{\sig} U_0^{\sig})
\Big\|_{S^2(E)} \\ & \le & \sum_{\sig \in \Gamma} \Degree
\|X_{\varepsilon}^{\sig}\|_{S_{\Degree}^2(S^2(E))}
\|\zeta^{\sig}(\omega_0) - U_0^{\sig}\|_{S_{\Degree}^2} + \Big\|
\sum_{\sig \in \Gamma} \Degree \mbox{tr}(X_{\varepsilon}^{\sig}
U_0^{\sig}) \Big\|_{S^2(E)} \\ & < & \,\ \delta \Delta_{\Gamma} +
\Big\| \sum_{\sig \in \Gamma} \Degree
\mbox{tr}(X_{\varepsilon}^{\sig} U_0^{\sig}) \Big\|_{S^2(E)}
\end{eqnarray*}
Taking $\varepsilon(\delta) = \delta / k_0(\delta)$, it is easy to
check that there exists $\gamma_1(\delta) > 0$ such that
\begin{equation} \label{Unitary}
\Big\| \sum_{\sig \in \Gamma} \Degree
\mbox{tr}(X_{\varepsilon(\delta)}^{\sig} U_0^{\sig})
\Big\|_{S^2(E)} \ge \Delta_{\Gamma} - \gamma_1(\delta)
\end{equation}
and where $\varepsilon(\delta), \gamma_1(\delta) \rightarrow 0^+$
as $\delta \rightarrow 0^+$. In particular, since for some other
election of the unitary matrices $U_0^{\sig}$ ($\sig \in \Gamma$)
we have the same value for $\lambda_{\sig} (\mathcal{U}_0(\sig,
\delta))$ (by the translation invariance of the Haar measure
$\lambda_{\sig}$), we obtain that $k_0(\delta)$ does not depend on
the chosen matrices $U_0^{\sig}$ ($\sig \in \Gamma$) and
(\ref{Unitary}) holds for any family of unitary matrices $U^{\sig}
\in U(\Degree)$ with $\sig \in \Gamma$. Now, given $A \in
\Lebesgue^1(\Sig)$ of norm 1, we use polar decomposition to write
$A^{\sig} = U_A^{\sig} |A^{\sig}|$ with $U_A^{\sig} \in
U(\Degree)$. Then, we have
\begin{eqnarray*}
\lefteqn{\Delta_{\Gamma} - \gamma_1(\delta)} \\ & \le & \Big\|
\sum_{\sig \in \Gamma} \Degree
\mbox{tr}(X_{\varepsilon(\delta)}^{\sig} U_A^{\sig} [I -
|A^{\sig}|]) \Big\|_{S^2(E)} + \Big\| \sum_{\sig \in \Gamma}
\Degree \mbox{tr}(X_{\varepsilon(\delta)}^{\sig} A^{\sig})
\Big\|_{S^2(E)} \\ & \le & \sum_{\sig \in \Gamma} \Degree
\|X_{\varepsilon(\delta)}^{\sig}\|_{S_{\Degree}^{\infty}(S^2(E))}
\|I - |A^{\sig}|\|_{S_{\Degree}^1} + \Big\| \sum_{\sig \in \Gamma}
\Degree \mbox{tr}(X_{\varepsilon(\delta)}^{\sig} A^{\sig})
\Big\|_{S^2(E)} \\ & \le & (1 + \xi(\delta))(\Delta_{\Gamma} - 1)
+ \Big\| \sum_{\sig \in \Gamma} \Degree
\mbox{tr}(X_{\varepsilon(\delta)}^{\sig} A^{\sig}) \Big\|_{S^2(E)}
\end{eqnarray*}
where the last inequality follows by Lemma \ref{Ultraproducts}.
Let us consider the subspace $F_{\Gamma}$ of $S^2(E)$ spanned by
the entries of $X_{\varepsilon(\delta)}^{\sig}$ where $\sig$ runs
over $\Gamma$. Then, the last inequality gives that the linear
isomorphism $\Lambda_{\Gamma}$, given by $$A \in
\Lebesgue^1(\Gamma) \longmapsto \sum_{\sig \in \Gamma} \Degree
\mbox{tr}(X_{\varepsilon(\delta)}^{\sig} A^{\sig}) \in
F_{\Gamma},$$ satisfies $\|\Lambda_{\Gamma}^{-1}\| \le 1 +
\gamma_2(\delta)$ for some $\gamma_2(\delta) > 0$ satisfying
$\gamma_2(\delta) \rightarrow 0^+$ as $\delta \rightarrow 0^+$.
Moreover, we have that $\|\Lambda_{\Gamma}\|_{cb} \le 1 +
\xi(\delta)$. Namely, given $A \in S^1(\Lebesgue^1(\Gamma))$, we
have
\begin{eqnarray*}
\Big\| \sum_{\sig \in \Gamma} \Degree \mbox{tr}
(X_{\varepsilon(\delta)}^{\sig} A^{\sig}) \Big\|_{S^1(S^2(E))} &
\le & \sum_{\sig \in \Gamma} \Degree
\|X_{\varepsilon(\delta)}^{\sig}\|_{S_{\Degree}^{\infty}(S^2(E))}
\|A^{\sig}\|_{S_{\Degree}^1(S^1)} \\ & \le & (1 + \xi(\delta)) \,\
\|A\|_{S^1(\Lebesgue^1(\Gamma))}.
\end{eqnarray*}
The first inequality follows by an inequality of Holder type, see
e.g. Lemma 3.3 of \cite{GP1}. That is, we have seen that
$\|\Lambda_{\Gamma}\|_{cb} = \|\Lambda_{\Gamma} \otimes I_{S^1}\|
\le 1 + \xi(\delta)$. Therefore $\|\Lambda_{\Gamma}\|_{cb}
\|\Lambda_{\Gamma}^{-1}\| \le 1 + \gamma_3(\delta)$ with
$\gamma_3(\delta) \rightarrow 0^+$ as $\delta \rightarrow 0^+$.
Therefore, taking $\delta \rightarrow 0^+$ we obtain that $E$
contains $\Lebesgue^1(\Gamma)$ $\lambda$-uniformly for all
$\lambda > 1$ as we wanted.
\end{proof}

\begin{remark}
\emph{As in Lemma \ref{Ultraproducts}, we could have used an
argument with ultrapro\-ducts to show that uniformly
non-$\Lebesgue^1(\Sig)$ operator spaces have $\Sig$-subtype. This
alternative proof is a bit shorter. However, for the shake of
clarity, we have preferred to give the more explicit argument used
in the proof of Theorem \ref{Equivalences}.}
\end{remark}

As it was pointed out in Remark \ref{Remark-p-Independence}, it is
very natural to wonder whether Definition \ref{Definitions} if
affected if we change $S^2(E)$ by $S^p(E)$ with $1 < p < \infty$.
The notion of $\B$-convexity should not depend on the election of
the exponent $p$ and, fortunately, this is the case.

\begin{corollary} \label{p-Convexity}
An operator space $E$ is $\B$-convex if and only if there exists a
finite subset $\Gamma$ of $\Sig$ such that
$$\frac{1}{\Delta_{\Gamma}} \inf_{B^{\sig} unitary}
\Big\|\sum_{\sig \in \Gamma} \Degree \mbox{tr}(A^{\sig} B^{\sig})
\Big\|_{S^p(E)} \le (1 - \delta) \,\ \max_{\sig \in
\Gamma}\|A^{\sig}\|_{S_{\Degree}^{\infty}(S^p(E))}$$ for some $1 <
p < \infty$ and any family $\{A^{\sig} \in M_{\Degree} \otimes
S^p(E)\}_{\sig \in \Gamma}$.
\end{corollary}

\begin{proof}
By Lemma 1.7 of \cite{P4} and Theorem \ref{Equivalences}, we know
that $E$ is $\B$-convex if and only if there exists a finite
subset $\Gamma$ of $\Sig$ such that $\|\mathrm{T}_2(\Gamma,
S^2(E))\| < \Delta_{\Gamma}^{1/2}$. On the other hand, given $1 <
p,q < \infty$, we claim that
\begin{equation} \label{p-Independence}
\|\mathrm{T}_2(\Gamma, S^p(E))\| < \Delta_{\Gamma}^{1/2}
\Longleftrightarrow \|\mathrm{T}_2(\Gamma, S^q(E))\| <
\Delta_{\Gamma}^{1/2}.
\end{equation}
By Holder inequality we have $\|\mathrm{T}_2(\Gamma, F)\| \le
\Delta_{\Gamma}^{1/2}$ for any operator space $F$. Then,
(\ref{p-Independence}) follows by complex interpolation with
$S^1(E)$ and $S^{\infty}(E)$. In particular, given $1 < p <
\infty$, we have that $E$ is $\B$-convex if and only if there
exists $\Gamma \subset \Sig$ finite such that
$\|\mathrm{T}_2(\Gamma, S^p(E))\| < \Delta_{\Gamma}^{1/2}$. But,
proceeding as in the proof of Theorem \ref{Equivalences}, we can
see that the desired inequality is equivalent to this last
condition.
\end{proof}

\begin{remark}
\emph{By similar arguments, if $1 < p < \infty$, we can also
replace $S^2(E)$ by $S^p(E)$ in the definition of uniformly
non-$\Lebesgue^1(\Sig)$ operator spaces.}
\end{remark}

\section{Tensor submultiplicativity}
\label{Section4}

As we pointed out in Remark \ref{Relevant-Parameters}, it seems
that the classical set of parameters and those $\Sig$ having
$\mathbf{d}_{\Sig}$ unbounded are the most relevant ones. In this
Section we shall see that, for these sets of parameters $\Sig$,
the notion of $\B$-convexity is stable under complete isomorphy
and the notion of containing $\Lebesgue^1(\Gamma)$
$\lambda$-uniformly does not depend on $\lambda > 1$. The
commutative analogs of these results are very well-known, see e.g.
\cite{G} or \cite{P1}. In order to prove these results, we need to
fix some notation. Let $(\Sig, \mathbf{d}_{\Sig})$ be a set of
parameters. Given two subsets $\Gamma^1$ and $\Gamma^2$ of $\Sig$,
we define their \emph{tensor product} as $$\Gamma^1 \otimes
\Gamma^2 = \{\sig_1 \otimes \sig_2: \,\ \sig_j \in \Gamma^j, \,\
j=1,2\} \qquad \mbox{where} \qquad d_{\sig_1 \otimes \sig_2} =
d_{\sig_1} d_{\sig_2}.$$ We say that $(\Sig, \mathbf{d}_{\Sig})$
is \emph{$\otimes$-closed} if, for any pair of subsets $\Gamma^1$
and $\Gamma^2$ of $\Sig$, there exists an injective mapping $j:
\Gamma^1 \otimes \Gamma^2 \rightarrow \Sig$ such that $d_{j(\sig_1
\otimes \sig_2)} = d_{\sig_1 \otimes \sig_2}$ for all $\sig_1
\otimes \sig_2 \in \Gamma^1 \otimes \Gamma^2$. Also, given $\Gamma
\subset \Sig$ finite, we define $$\mathrm{N}_{\Gamma}(E) =
\frac{1}{\sqrt{\Delta_{\Gamma}}} \,\
\|\mathrm{T}_2(\Gamma,E)\|_{cb}.$$

\begin{lemma} \label{Submultiplicativity}
If $(\Sig, \mathbf{d}_{\Sig})$ is $\otimes$-closed, then
$\mathrm{N}_{\Gamma^1 \otimes \Gamma^2}(E) \le
\mathrm{N}_{\Gamma^1}(E) \,\ \mathrm{N}_{\Gamma^2}(E)$ for any
pair of finite subsets  $\Gamma^1$ and $\Gamma^2$ of $\Sig$.
\end{lemma}

\begin{proof}
Let us consider a family $A = \{A^{\sig_1 \otimes \sig_2} \in
M_{d_{\sig_1}d_{\sig_2}} \otimes S^2(E): \,\ \sig_1 \in \Gamma^1,
\sig_2 \in \Gamma^2\}$. Let $A^{\sig_1 \otimes \sig_2}(\omega) =
\zeta^{j(\sig_1 \otimes \sig_2)}(\omega) A^{\sig_1 \otimes
\sig_2}$ for $\omega \in \Omega$. Since $\zeta^{j(\sig_1 \otimes
\sig_2)}$ is uniformly distributed on the unitary group
$U(d_{\sig_1}d_{\sig_2})$ and $\zeta^{\sig_1}(\omega_1) \otimes
\zeta^{\sig_2}(\omega_2)$ is unitary, we have
\begin{eqnarray*}
\lefteqn{\int_{\Omega} \Big\| \sum_{\sig_j \in \Gamma^j}
d_{\sig_1} d_{\sig_2} \mbox{tr}(A^{\sig_1 \otimes \sig_2}(\omega))
\Big\|_{S^2(E)}^2 d \mu(\omega)} \\ & = & \int_{\Omega} \Big\|
\sum_{\sig_j \in \Gamma^j} d_{\sig_1} d_{\sig_2}
\mbox{tr}(A^{\sig_1 \otimes \sig_2}(\omega)
(\zeta^{\sig_1}(\omega_1) \otimes \zeta^{\sig_2}(\omega_2)))
\Big\|_{S^2(E)}^2 d\mu(\omega)
\end{eqnarray*}
for all $\omega_1, \omega_2 \in \Omega$. Therefore, if we write
$$X^{\sig_1} (\omega, \omega_2) = \sum_{\sig_2 \in \Gamma^2}
d_{\sig_2} \mbox{tr} (A^{\sig_1 \otimes \sig_2}(\omega)
\zeta^{\sig_2}(\omega_2)) \in M_{d_{\sig_1}} \otimes S^2(E),$$ we
obtain the following estimate
\begin{eqnarray*}
\lefteqn{\int_{\Omega} \Big\| \sum_{\sig_j \in \Gamma^j}
d_{\sig_1} d_{\sig_2} \mbox{tr}(A^{\sig_1 \otimes \sig_2}(\omega))
\Big\|_{S^2(E)}^2 d \mu(\omega)} \\ & = & \! \! \int_{\Omega^3}
\Big\| \sum_{\sig_j \in \Gamma^j} d_{\sig_1} d_{\sig_2}
\mbox{tr}(A^{\sig_1 \otimes \sig_2}(\omega)
(\zeta^{\sig_1}(\omega_1) \otimes \zeta^{\sig_2}(\omega_2)))
\Big\|_{S^2(E)}^2 d\mu(\omega_1) d\mu(\omega_2) d\mu(\omega) \\ &
= & \! \! \int_{\Omega^3} \Big\| \sum_{\sig_1 \in \Gamma^1}
d_{\sig_1} \mbox{tr} (X^{\sig_1}(\omega, \omega_2)
\zeta^{\sig_1}(\omega_1)) \Big\|_{S^2(E)}^2 d\mu(\omega_1)
d\mu(\omega_2) d\mu(\omega) \qquad \qquad {}
\end{eqnarray*}
\begin{eqnarray*}
& \le & \mathrm{N}_{\Gamma^1}(E)^2 \, \Delta_{\Gamma^1} \,
\sum_{\sig_1 \in \Gamma^1} d_{\sig_1} \int_{\Omega \times \Omega}
\|X^{\sig_1}(\omega, \omega_2)\|_{S_{d_{\sig_1}}^2(S^2(E))}^2
d\mu(\omega_2) d\mu(\omega) \\ & \le & \mathrm{N}_{\Gamma^1}(E)^2
\, \mathrm{N}_{\Gamma^2}(E)^2 \, \Delta_{\Gamma^1}
\Delta_{\Gamma^2} \, \sum_{\sig_j \in \Gamma^j} d_{\sig_1}
d_{\sig_2} \int_{\Omega} \|A^{\sig_1 \otimes
\sig_2}(\omega)\|_{S_{d_{\sig_1}d_{\sig_2}}^2(S^2(E))}^2
d\mu(\omega) \\ & = & \mathrm{N}_{\Gamma^1}(E)^2 \,
\mathrm{N}_{\Gamma^2}(E)^2 \, \Delta_{\Gamma^1} \Delta_{\Gamma^2}
\, \sum_{\sig_j \in \Gamma^j} d_{\sig_1} d_{\sig_2} \|A^{\sig_1
\otimes \sig_2}\|_{S_{d_{\sig_1}d_{\sig_2}}^2(S^2(E))}^2
\end{eqnarray*}
The last equality follows by the unitarity of $\zeta^{j(\sig_1
\otimes \sig_2)}(\omega)$. Now, since $\Delta_{\Gamma^1 \otimes
\Gamma^2} = \Delta_{\Gamma^1} \Delta_{\Gamma^2}$ we obtain the
desired inequality. This completes the proof.
\end{proof}

\begin{proposition} \label{Stability}
If $(\Sig, \mathbf{d}_{\Sig})$ is $\otimes$-closed, then
$\B$-convexity is stable under complete isomorphism.
\end{proposition}

\begin{proof}
Let us assume that $E$ is $\B$-convex and let $F$ be an operator
space completely isomorphic to $E$. By Theorem \ref{Equivalences},
we know that there exists some $\Gamma_0 \subset \Sig$ finite such
that $\mathrm{N}_{\Gamma_0}(E) < 1$ and it suffices to see that
there exists some $\Gamma \subset \Sig$ finite such that
$\mathrm{N}_{\Gamma}(F) < 1$. However, if $d_{cb}(E,F)$ stands for
the $cb$-distance between $E$ and $F$, we have
$$\mathrm{N}_{\Gamma_0^{\otimes n}}(F) \le d_{cb}(E,F) \,\
\mathrm{N}_{\Gamma_0^{\otimes n}}(E)$$ where $n$ is any positive
integer and $\Gamma_0^{\otimes n} = \Gamma_0 \otimes \Gamma_0
\otimes \cdots \otimes \Gamma_0$ with $n$ factors. Now, by Lemma
\ref{Submultiplicativity}, we know that
$\mathrm{N}_{\Gamma_0^{\otimes n}}(E) \le
\mathrm{N}_{\Gamma_0}(E)^n$. Finally, we are done by taking
$\Gamma = \Gamma_0^{\otimes n}$ with $n$ large enough. This
completes the proof.
\end{proof}

\begin{remark} \label{Remark-Stability}
\emph{Proposition \ref{Stability} obviously holds for the
classical set of parameters. On the other hand, if we consider the
set of parameters given by $\Sig = \N \times \N$ and
$d_{\sig_{jk}} = 2^k$, then $\sig_{11}$ generates the whole set of
parameters by taking tensor powers of it. This set satisfies that
$\mathbf{d}_{\Sig}$ is unbounded and is again $\otimes$-closed.
Therefore, Proposition \ref{Stability} also holds for it.
Moreover, in Section \ref{Section6} we shall prove that the notion
of $\B$-convexity does not depend on $\Sig$ whenever we work with
sets of parameters with $\mathbf{d}_{\Sig}$ unbounded. In
particular, we have seen that Proposition \ref{Stability} holds
for any set of parameters with $\mathbf{d}_{\Sig}$ unbounded.}
\end{remark}

\begin{proposition} \label{Lambda-Independence}
Let us suppose that $(\Sig, \mathbf{d}_{\Sig})$ is
$\otimes$-closed and let $E$ be an operator space containing
$\Lebesgue^1(\Gamma)$ $\lambda$-uniformly for some $\lambda > 1$.
Then, for all $\tau > 1$, $E$ contains $\Lebesgue^1(\Gamma)$
$\tau$-uniformly.
\end{proposition}

\begin{proof}
We have already seen in the proof of Theorem \ref{Equivalences}
that, if $E$ contains $\Lebesgue^1 (\Gamma)$ $\lambda$-uniformly,
we obtain $$\frac{1}{\Delta_{\Gamma}} \inf_{B^{\sig} unitary}
\Big\| \sum_{\sig \in \Gamma} \Degree \mbox{tr}(X^{\sig} B^{\sig})
\Big\|_{S^2(E)} \ge \frac{1}{\lambda} \max_{\sig \in \Gamma}
\|X^{\sig}\|_{S_{\Degree}^{\infty}(S^2(E))}$$ for all $\Gamma
\subset \Sig$ finite and certain family of matrices $X^{\sig} \in
M_{\Degree} \otimes S^2(E)$, with $\sig \in \Gamma$. On the other
hand, from this inequality it is not difficult to see using Lemma
\ref{cb-distance} that $\mathrm{N}_{\Gamma}(E) \ge 1 / \lambda$
for all $\Gamma \subset \Sig$ finite. Now, if $E$ is $\B$-convex,
we know that there exist some finite subset $\Gamma_0$ of $\Sig$
such that $\mathrm{N}_{\Gamma_0}(E) < 1$. Therefore, by Lemma
\ref{Submultiplicativity}, we would have $$\frac{1}{\lambda} \le
\mathrm{N}_{\Gamma_0^{\otimes n}}(E) \le
\mathrm{N}_{\Gamma_0}(E)^n \rightarrow 0^+ \qquad \mbox{as} \ \ n
\rightarrow \infty.$$ This gives that $E$ is not $\B$-convex. But,
by Theorem \ref{Equivalences}, this is equivalent to say that $E$
contains $\Lebesgue^1(\Gamma)$ $\tau$-uniformly for all $\tau
> 1$. This completes the proof.
\end{proof}

\begin{remark} \label{Remark-Lambda}
\emph{Proposition \ref{Lambda-Independence} obviously holds for
the classical set of parameters. Moreover, by similar arguments
that the ones given in Remark \ref{Remark-Stability}, Proposition
\ref{Lambda-Independence} also holds for any set of parameters
$(\Sig, \mathbf{d}_{\Sig})$ with $\mathbf{d}_{\Sig}$ unbounded.}
\end{remark}

\section{Non-trivial $\Sig$-type and $\B$-convexity}
\label{Section5}

In this Section, $\Sig_0$ will stand for the classical set of
parameters. We begin by showing that any operator space having
non-trivial $\Sig$-type is $\B$-convex. However the most
interesting point is that, in contrast with the classical theory,
the converse is false. We shall provide examples of
$\mathrm{B}_{\Sig_0}$-convex operator spaces failing to have
$\Sig_0$-type $p$ for any $1 < p \le 2$. This is a very important
difference between the commutative and non-commutative contexts.
Namely, it turns out that we can not expect to obtain an operator
space version of the Maurey-Pisier theorem. We recall that this
result asserts that, for any infinite-dimensional Banach space
$B$, the supremum of those $p \in [1,2]$ for which $B$ has type
$p$ coincides with the minimum of those $1 \le q \le 2$ for which
$B$ contains $l_n^q$ uniformly, see \cite{MaP} for more details.

The examples we are giving are the well-known row and column
operator spaces $R$ and $C$, see e.g. \cite{P3} for the definition
of these spaces. This is even more surprising since $R$ and $C$
are Hilbertian. Moreover, we shall provide some other examples of
Hilbertian operator spaces having sharp $\Sig_0$-type $p$ for any
$1 < p \le 2$, which are obviously $\mathrm{B}_{\Sig_0}$-convex.
Finally, we shall use a result of Pisier to show that $\min l^2$
and $\max l^2$ are Hilbertian operator spaces failing the
$\mathrm{B}_{\Sig_0}$-convexity.

\begin{proposition} \label{Weak-Type}
If $E$ has non-trivial $\Sig$-type, then $E$ is $\B$-convex.
\end{proposition}

\begin{proof}
Let $\Gamma$ be any finite subset of $\Sig$ and let us suppose
that $E$ has $\Sig$-type $p$ for some $1 < p \le 2$. By the
operator space version of the classical Minkowski inequality (see
e.g. \cite{GP1}), we have that the natural mapping
$$\Lebesgue_{S^2(E)}^p(\Gamma) \longrightarrow
S^2(\Lebesgue_E^p(\Gamma))$$ is completely contractive. In
particular, if $\mathcal{K}_p(E, \Sig)$ stands for the $cb$-norm
of $\mathrm{T}_p(E)$, we have $\|\mathrm{T}_p(\Gamma, S^2(E))\|
\le \mathcal{K}_p(E, \Sig).$ Hence,
$$\|\mathrm{T}_2(\Gamma,E)\|_{cb} = \|\mathrm{T}_2(\Gamma, S^2(E)
)\| \le \|\mathrm{T}_p(\Gamma, S^2(E))\| \,\ \Delta_{\Gamma}^{1/p
- 1/2} \le \mathcal{K}_p(E, \Sig) \,\ \Delta_{\Gamma}^{1/p -
1/2},$$ where the first inequality follows easily from Lemma
\ref{cb-distance}. Therefore, taking $\Gamma$ large enough so that
$\mathcal{K}_p(E, \Sig) < \Delta_{\Gamma}^{1 - 1/p}$, we conclude
that $E$ has $\Sig$-subtype. Thus, $E$ is $\B$-convex by Theorem
\ref{Equivalences}. This completes the proof.
\end{proof}

The following result has its origins in an unpublished result of
Magdalena Musat which asserts that $S^2(R)$ and $S^2(C)$ are
superreflexive Banach spaces. After some conversations, initiated
by Marius Junge and Gilles Pisier, Timur Oikhberg found a
surprisingly simple proof of this fact. The next Theorem is based
on the techniques employed there. First, we fix some notation. As
usual, given $0 < \theta < 1$, we shall denote by $R(\theta)$ the
complex interpolation operator space $(R,C)_{\theta}$.
Analogously, $C(\theta)$ stands for $(C,R)_{\theta} = R(1 -
\theta)$. By convention, we also set that $R(0) = C(1) = R$ and
$C(0) = R(1) = C$.

\begin{theorem} \label{Counterexample}
Let $1 \le p \le 2$, then $R(1/p)$ and $C(1/p)$ are
$\mathrm{B}_{\Sig_0}$-convex Hilbertian operator spaces having
$\Sig_0$-type $p$. Moreover, if $1 \le p < 2$, then $R(1/p)$ and
$C(1/p)$ do not have $\Sig_0$-type $q$ for any $p < q \le 2$.
\end{theorem}

\begin{proof}
Let us suppose that $R(1/p)$ and $C(1/p)$ are not
$\mathrm{B}_{\Sig_0}$-convex. Then, by Theorem \ref{Equivalences}
and Remark \ref{Uniformly-Classical}, the spaces $S^2(R(1/p))$ and
$S^2(C(1/p))$ should not have type $> 1$ in the Banach space
sense. However, we claim that both spaces have type $4/3$.
Therefore, $R(1/p)$ and $C(1/p)$ are $\mathrm{B}_{\Sig_0}$-convex.
To prove our claim we recall that, by Theorem 1.1 of \cite{P4}, we
have $S^2(R) = R(1/2) \otimes_h R \otimes_h R(1/2)$ and $S^2(C) =
C(1/2) \otimes_h C \otimes_h C(1/2)$ completely isometrically.
Now, since the Haagerup tensor product commutes with the complex
interpolation functor, we can write
\begin{equation} \label{Complex-Interpolation}
\begin{array}{c} S^2(R) = (R \otimes_h R \otimes_h R, C \otimes_h R
\otimes_h C)_{1/2} \\ S^2(C) = (C \otimes_h C \otimes_h C, R
\otimes_h C \otimes_h R)_{1/2} \end{array}
\end{equation}
completely isometrically. But, as Banach spaces, $R \otimes_h R
\otimes_h R$ and $C \otimes_h C \otimes_h C$ are isometrically
isomorphic to a Hilbert space. In particular,
(\ref{Complex-Interpolation}) gives that $S^2(R)$ and $S^2(C)$
have type $4/3$ in the Banach space sense. Hence, by complex
interpolation, the same happens with $S^2(R(1/p))$ and
$S^2(C(1/p))$.

On the other hand, by the reiteration theorem for the complex
interpolation method, we have $R(1/p) = (R(1/2), C)_{\frac{2}{p} -
1}$ and $C(1/p) = (C(1/2), R)_{\frac{2}{p} - 1}$. But $R(1/2) =
C(1/2)$ is an OH operator space, see \cite{P3}. Therefore, it is
easy to check that $R(1/2)$ has $\Sig_0$-type $2$, see \cite{GP2}
for more on this topic. In particular, by complex interpolation,
we get that $R(1/p)$ and $C(1/p)$ have at least $\Sig_0$-type $q$
where $$\frac{1}{q} = \frac{1 - (2/p - 1)}{2} + \frac{2/p - 1}{1}
= \frac{1}{p}.$$ That is, $R(1/p)$ and $C(1/p)$ have $\Sig_0$-type
$p$. Now let $p < q \le 2$, we want to see that $R(1/p)$ and
$C(1/p)$ do not have $\Sig_0$-type $q$. Following the notation
introduced in Remark \ref{Independence-Type}, we obviously have
that $$\|\mathrm{T}_q(\Sig_0, E)\|_{cb} \ge
\|\mathrm{T}_q^q(\Sig_0, E)\|_{cb} = \|\mathrm{T}_q^q(\Sig_0,
S^q(E))\|$$ for any operator space $E$. Hence, we just need to
check that $\mathrm{T}_q^q(\Sig_0, S^q(E))$ is not bounded if $E =
R(1/p)$ or $E = C(1/p)$. But, by the Khintchine-Kahane
inequalities, this is equivalent to saying that $S^q(R(1/p))$ and
$S^q(C(1/p))$ do not have type $q$ in the Banach space sense for
any $p < q \le 2$. Let us recall that
\begin{eqnarray*}
S^q(R(1/p)) & = & (S^{\infty}(R(1/p)), S^1(R(1/p)))_{1/q} \\ & = &
(C \otimes_h R(1/p) \otimes_h R, R \otimes_h R(1/p) \otimes_h
C)_{1/q} \\ & = & C(1/q) \otimes_h R(1/p) \otimes _h R(1/q).
\end{eqnarray*}
Analogously, we obtain $S^q(C(1/p)) = C(1/q) \otimes_h C(1/p)
\otimes _h R(1/q)$. Let us consider the subspace of $S^q(R(1/p))$
corresponding to $C(1/q) \otimes_h R(1/p)$. Then we can write
\begin{eqnarray*}
C(1/q) \otimes_h R(1/p) & = & (C \otimes_h R(1/p), R \otimes_h
R(1/p))_{1/q} \\ & = & ((C \otimes_h R, C \otimes_h C)_{1/p}, (R
\otimes_h R, R \otimes_h C)_{1/p})_{1/q}
\end{eqnarray*}
completely isometrically. But, as Banach spaces, $R \otimes_h R$
and $C \otimes_h C$ are isometrically isomorphic to $S^2$.
Therefore, we have $$C(1/q) \otimes_h R(1/p) = ((S^{\infty},
S^2)_{1/p}, (S^2, S^1)_{1/p})_{1/q} = (S^{2p}, S^{2p/p+1})_{1/q} =
S^{2pq/p+q}$$ isometrically. Finally, since $2pq < q (p+q)$
whenever $p < q$, we have that $C(1/q) \otimes_h R(1/p)$ can not
have type $q$ in the Banach space sense. Consequently, the same
happens for $S^q(R(1/p))$. A similar argument gives that
$S^q(C(1/p))$ can not have type $q$ in the Banach space sense.
This completes the proof.
\end{proof}

\begin{remark}
\emph{In particular, by Theorem \ref{Counterexample}, the row and
column operator spaces are examples of
$\mathrm{B}_{\Sig_0}$-convex Hilbertian operator spaces failing to
have non-trivial $\Sig_0$-type.}
\end{remark}

\begin{remark}
\emph{Although we give more details in Section \ref{Section6}, it
is a simple consequence of Theorem \ref{Equivalences} that
$\mathrm{B}_{\Sig_0}$-convexity is the strongest condition among
the possible sets of parameters we are working with. That is, a
$\mathrm{B}_{\Sig_0}$-convex operator space is automatically
$\B$-convex for any other set of parameters $\Sig$. In particular,
the examples treated in Theorem \ref{Counterexample} are
$\B$-convex. Moreover, the given argument to see that $R(1/p)$ and
$C(1/p)$ have $\Sig_0$-type $p$ for any $1 \le p \le 2$ remains
valid for any other set of parameters $\Sig$.}
\end{remark}

Once we have found examples of $\mathrm{B}_{\Sig_0}$-convex
Hilbertian operator spaces having sharp $\Sig_0$-type $p$ for any
$1 \le p \le 2$, we now show that $\min l^2$ and $\max l^2$ are
Hilbertian operator spaces failing to be
$\mathrm{B}_{\Sig_0}$-convex. This is based on Example 4.2 of
\cite{P4}, where Pisier makes the following construction. Let
$M_2$ be the algebra of $2 \times 2$ complex-valued matrices
equipped with its normalized trace $t$ and let us set $(A_k, t_k)
= (M_2, t)$ for any $k \ge 1$. Then we consider the so-called
hyperfinite $\mathrm{II}_1$ factor $$(\mathcal{M}, \tau) =
\bigotimes_{k=1}^{\infty} (A_k, t_k).$$ Let $\mathcal{M}_n$ stand
for the subalgebra of $\mathcal{M}$ corresponding to $A_1 \otimes
\cdots \otimes A_n$. Let us consider the element of $\mathcal{M}
\otimes \min l^2$ given by $d_n = V_n \otimes e_n$ where we write
$(e_i)$ for the canonical basis of $l^2$ and $(V_n)$ is a sequence
in $\mathcal{M}$ satisfying $V_n \in \mathcal{M}_n$ for all $n \ge
1$, $E^{\mathcal{M}_n}(V_{n+1}) = 0$ and the canonical
anticommutation relations $$V_i V_j^{\star} + V_j^{\star} V_i =
\delta_{ij} I \qquad \mbox{and} \qquad V_i V_j + V_j V_i = 0.$$
Then, from these properties, Pisier shows that for all finite
sequence of scalars $(\alpha_k)$ with $1 \le k \le n$ and for all
$1 \le p \le \infty$ we have
\begin{equation} \label{min-l2}
\frac{1}{2} \,\ \sup_{1 \le k \le n} |\alpha_k| \le \Big\| \sum_{k
= 1}^n \alpha_k d_k \Big\|_{L^p(\mathcal{M}; \min l^2)} \le
\sup_{1 \le k \le n} |\alpha_k|
\end{equation}
where $L^p(\mathcal{M}; \min l^2)$ denotes the non-commutative
$L^p$ space defined in $(\mathcal{M}, \tau)$ and with values in
$\min l^2$. In particular, inequalities (\ref{min-l2}) tell us
that the Banach-Mazur distance between some subspace of
$L^p(\mathcal{M}_n; \min l^2)$ and $l_n^{\infty}$ is bounded above
by $2$ for all $n \ge 1$ and all $1 \le p \le \infty$. Then,
recalling the natural embedding of $l_n^1$ into $l_{2^n}^{\infty}$
and taking $p = 2$ in (\ref{min-l2}), it is easy to see that
$S^2(\min l^2)$ contains $l_n^1$ uniformly in the Banach space
sense. Therefore, by Remark \ref{Uniformly-Classical} and Theorem
\ref{Equivalences}, we have that $\min l^2$ is not
$\mathrm{B}_{\Sig_0}$-convex. By a duality argument, the same
happens for $\max l^2$.

\section{On the independence with respect to $\Sig$}
\label{Section6}

In this last Section we study the dependence of the notion of
$\B$-convexity with respect to the set of parameters $\Sig$. We
begin by showing the independence with respect to $\Sig$ when we
work with sets of parameters satisfying that $\mathbf{d}_{\Sig}$
is an unbounded family. After that, we shall give two interesting
equivalent formulations of the possible independence of
$\B$-convexity with respect to any set $\Sig$.

\begin{proposition} \label{Sig-Independence}
Let us consider two sets of parameters $(\Sig_1,
\mathbf{d}_{\Sig_1})$ and $(\Sig_2, \mathbf{d}_{\Sig_2})$ with
$\mathbf{d}_{\Sig_1}$ and $\mathbf{d}_{\Sig_2}$ unbounded. Let $E$
be an operator space, then $E$ is $\mathrm{B}_{\Sig_1}$-convex if
and only if $E$ is $\mathrm{B}_{\Sig_2}$-convex.
\end{proposition}

\begin{proof}
By Theorem \ref{Equivalences}, we just need to check that $E$
contains $\Lebesgue^1(\Gamma)$ $\lambda$-uniformly ($\Gamma
\subset \Sig_1$) for all $\lambda > 1$ if and only if $E$ contains
$\Lebesgue^1(\Gamma)$ $\lambda$-uniformly ($\Gamma \subset
\Sig_2$) for all $\lambda > 1$. In particular, it suffices to see
that $\Lebesgue^1(\Sig_1)$ contains $\Lebesgue^1(\Gamma)$
$\lambda$-uniformly ($\Gamma \subset \Sig_2$) for all $\lambda >
1$ and that $\Lebesgue^1(\Sig_2)$ contains $\Lebesgue^1(\Gamma)$
$\lambda$-uniformly ($\Gamma \subset \Sig_1$) for all $\lambda >
1$. But this follows from the unboundedness of
$\mathbf{d}_{\Sig_1}$ and $\mathbf{d}_{\Sig_2}$. Namely, given
$\Gamma \subset \Sig_2$ finite, we know that there exists $\Lambda
\subset \Sig_1$ and a bijection $\tau: \Gamma \rightarrow \Lambda$
such that $\Degree \le d_{\tau(\sig)}$ for all $\sig \in \Gamma$.
In particular we can consider the linear mapping
$\mathrm{S}_{\Gamma}: \Lebesgue^1(\Gamma) \rightarrow
\Lebesgue^1(\Lambda)$ given by
$$\mathrm{S}_{\Gamma}(A)_{ij}^{\tau(\sig)} =
\frac{\Degree}{d_{\tau(\sig)}} \left\{ \begin{array}{ll}
A_{ij}^{\sig} & \mbox{if} \ \ 1 \le i,j \le \Degree \\ 0 &
\mbox{otherwise}. \end{array} \right.$$ Given the fact that
$\mathrm{S}_{\Gamma}$ is a complete isometry,
$\Lebesgue^1(\Sig_1)$ contains $\Lebesgue^1(\Gamma)$ $1$-uniformly
(where $\Gamma \subset \Sig_2$). Similarly, using the
unboundedness of $\mathbf{d}_{\Sig_2}$, we can see that
$\Lebesgue^1(\Sig_2)$ contains $\Lebesgue^1(\Gamma)$ $1$-uniformly
($\Gamma \subset \Sig_1$). This completes the proof.
\end{proof}

\begin{remark} \label{Sig-Order}
\emph{Given two sets of parameters $(\Sig_1, \mathbf{d}_{\Sig_1})$
and $(\Sig_2, \mathbf{d}_{\Sig_2})$, we shall say that $\Sig_1 \le
\Sig_2$ if there exists an injective mapping $j: \Sig_1
\rightarrow \Sig_2$ such that $\Degree \le d_{j(\sig)}$ for all
$\sig \in \Sig_1$. Then, by similar arguments to those used in
Proposition \ref{Sig-Independence}, it is easy to see that
$\mathrm{B}_{\Sig_1}$-convexity is stronger than
$\mathrm{B}_{\Sig_2}$-convexity whenever $\Sig_1 \le \Sig_2$. In
particular, if $\Sig_0$ stands for the classical set of
parameters, a $\mathrm{B}_{\Sig_0}$-convex operator space is
automatically $\B$-convex for any other set of parameters $\Sig$.}
\end{remark}

Proposition \ref{Sig-Independence} is only a little step in order
to see the independence of the notion of $\B$-convexity with
respect to $\Sig$. However, the general case seems to be more
complicated. Now, we give two different conditions which could be
useful to decide whether or not the independence with respect to
$\Sig$ holds.

\begin{itemize}
\item[(A)] \textbf{On the notion of $\K$-convexity.} In order to
introduce $\K$-convexity, we need to define the quantized version
of the Gauss system. It was already defined by Marcus and Pisier
in \cite{MP}. More in connection with the present paper, this
system was also treated in the operator space version of
Kwapie\'n's theorem, see \cite{GP2}. Given a set of parameters
$(\Sig, \mathbf{d}_{\Sig})$, we consider a family of independent
standard complex-valued gaussian random variables
$\{\gamma_{ij}^{\sig}: \Omega \rightarrow \C: \sig \in \Sig, \,\ 1
\le i,j \le \Degree\}$ indexed by $\Sig$ and $\mathbf{d}_{\Sig}$.
Then, if we construct the random matrices $$\gamma^{\sig} =
\frac{1}{\sqrt{\Degree}} \,\ \Big( \,\ \gamma_{ij}^{\sig} \,\
\Big),$$ we obtain the \emph{quantized Gaussian system of
parameters} $(\Sig, \mathbf{d}_{\Sig})$. On the other hand, given
an operator space $E$ and $f \in L_E^2(\Omega)$, we can consider
the Fourier coefficients of $f$ with respect to this system
$$\widehat{f}(\sig) = \int_{\Omega} f(\omega)
\gamma^{\sig}(\omega)^{\star} d \mu(\omega) \in M_{\Degree}
\otimes E.$$ We shall say that an operator space $E$ is
\emph{$\K$-convex} if the gaussian projection defined as follows
$$f \in L_E^2(\Omega) \longmapsto \sum_{\sig \in \Sig} \Degree
\mbox{tr}(\widehat{f}(\sig) \gamma^{\sig}) \in L_E^2(\Omega)$$ is
a completely bounded mapping. However, it is obvious that
$$\sum_{\sig \in \Sig} \Degree \mbox{tr}(\widehat{f}(\sig)
\gamma^{\sig}) = \sum_{\sig \in \Sig} \sum_{i,j = 1}^{\Degree}
\int_{\Omega} f(\omega) \overline{\gamma_{ij}^{\sig}(\omega)} d
\mu(\omega) \,\ \gamma_{ij}^{\sig}$$ where $\gamma_{ij}^{\sig}$
are independent complex-valued gaussian random variables. Hence,
it turns out that the notion of $\K$-convexity does not depend on
the set of parameters $\Sig$. Moreover, if $\Sig_0$ stands for the
classical set of parameters, then any operator space $E$ satisfies
\begin{eqnarray*}
E \ \ \mbox{$\K$-convex} & \Longleftrightarrow & E \ \
\mbox{$\mathrm{K}_{\Sig_0}$-convex} \\ & \Longleftrightarrow &
S^2(E) \ \ \mbox{$\mathrm{K}$-convex as a Banach space} \\ &
\Longleftrightarrow & S^2(E) \ \ \mbox{$\mathrm{B}$-convex as a
Banach space} \\ & \Longleftrightarrow & E \ \
\mbox{$\mathrm{B}_{\Sig_0}$-convex}.
\end{eqnarray*}
Therefore, it follows that the notion of $\B$-convexity does not
depend on the set of parameters $\Sig$ if and only if
$\B$-convexity and $\K$-convexity are equivalent notions. In
particular, it provides a possible approach to check this
independence. That is, the problem is to generalize to the
operator space setting Pisier's theorem which shows that
$\mathrm{K}$-convex and $\mathrm{B}$-convex Banach spaces are the
same, see \cite{P2} for more details on this topic.

\item[(B)] \textbf{On how $S_n^1$ embeds in $S^2(l^1)$.}
Let us suppose that we are given a subspace $F_n$ of
$S^2(l_{n^2}^1)$ and a linear isomorphism $\Phi_n: S_n^1
\rightarrow F_n$ for each $n \ge 1$. Let us denote by $F_n^{\min}$
the subspace $F_n$ equipped with the operator space structure
inherited as a subspace of $S^2(\min l_{n^2}^1)$. Then, if we
write $\Psi_n: F_n^{\min} \rightarrow S_n^1$ for the inverse of
$\Phi_n$ with the modified operator space structure on $F_n$, we
claim that the condition
\begin{equation} \label{Control-Min}
\|\Phi_n\|_{cb} \|\Psi_n\| \le k \qquad \mbox{for all $n \ge 1$
and some constant $k > 1$}
\end{equation}
implies the $\Sig$-independence of $\B$-convexity. Namely, let $E$
be an operator space which contains $l^1_n$ uniformly in the sense
of Definition \ref{Definitions}. That is, for all $n \ge 1$ there
exists a subspace $K_n$ of $S^2(E)$ and a linear isomorphism
$\Lambda_n: l_n^1 \rightarrow K_n$ such that $\|\Lambda_n\|_{cb}
\|\Lambda_n^{-1}\| \le \lambda$. Now, if we consider the linear
isomorphism $$\mathrm{T}_n = (I_{S^2} \otimes \Lambda_{n^2}) \circ
\Phi_n: S_n^1 \longrightarrow H_n \subset S^2(E),$$ the inverse
operator factors as follows $$H_n \longrightarrow
F_n^{\mbox{\scriptsize{min}}} \longrightarrow S_n^1$$ via de
composition $\mathrm{T}_n^{-1} = \Psi_n \circ (I_{S^2} \otimes
\Lambda_{n^2}^{-1})$. In summary, using the well-known properties
of the minimal operator space structure, the following estimate
follows from condition (\ref{Control-Min}) $$\|\mathrm{T}_n\|_{cb}
\|\mathrm{T}_n^{-1}\| \le \|\Lambda_{n^2}\|_{cb} \|\Phi_n\|_{cb}
\|\Psi_n\| \|\Lambda_{n^2}^{-1}\| \le k \lambda.$$ This gives that
$E$ contains $\Lebesgue^1(\Gamma)$ uniformly for any set of
parameters $(\Sig, \mathbf{d}_{\Sig})$. Therefore, we have proved
that any $\B$-convex operator space is
$\mathrm{B}_{\Sig_0}$-convex. In particular, the
$\Sig$-independence of $\B$-convexity follows from Remark
\ref{Sig-Order}.
\end{itemize}

\begin{remark}
\emph{As in the Banach space case, the given definition of
$\K$-convexity should not change if we consider the gaussian
projection on $L_E^p(\Omega)$ instead of $L_E^2(\Omega)$ for any
$1 < p < \infty$. Fortunately, this is the case. However, this
time the proof can not be supported by the Khintchine-Kahane
inequalities as it was explained in Remark
\ref{Khintchine-Kahane}. Nevertheless, the argument is simple.
Namely, if we refer to this a priori new notion as
$\K^p$-convexity, then an operator space $E$ is $\K^p$-convex if
and only if it is $\mathrm{K}_{\Sig_0}^p$-convex. But this last
condition means that $S^p(E)$ is $\mathrm{K}$-convex as a Banach
space. Now, since $\mathrm{K}$-convex and $\mathrm{B}$-convex
Banach spaces are the same, we conclude that $E$ is $\K^p$-convex
if and only if $S^p(E)$ is $\mathrm{B}$-convex as a Banach space.
But, by Corollary \ref{p-Convexity}, we know that this is
equivalent to say that $E$ is $\mathrm{B}_{\Sig_0}$-convex.
Finally, since $\K$-convexity is equivalent to
$\mathrm{B}_{\Sig_0}$-convexity, we are done.}
\end{remark}

\begin{remark}
\emph{Condition (\ref{Control-Min}) holds and so Pisier's theorem
on the equivalence bet\-ween B-convex and K-convex spaces remains
valid in the non-commutative setting. The proof of this fact,
which will appear as part of a joint work with Marius Junge
\cite{JP}, became clear after this paper was submitted for
publication.}
\end{remark}

\

\textbf{Acknowledgements.} The author wishes to thank Marius
Junge, Magdalena Musat and Gilles Pisier for many illuminating
conversations and some suggestions concerning the content of this
paper. Research supported in part by the European Commission via
the IHP Network \lq Harp\rq ${}$ and by MCYT Spain via the Project
BFM 2001/0189.

\bibliographystyle{amsplain}

\end{document}